\documentclass[11pt]{article}

\usepackage[T1]{fontenc}

\usepackage{amsmath,amsthm,amssymb}
\usepackage{graphicx}
\usepackage[all]{xy}
\usepackage{psfrag}

\usepackage{color}
\usepackage{epsfig}
\usepackage{latexsym}

\newtheorem{thm}{Theorem}[section]
\newtheorem*{thmintro}{Theorem}

\newtheorem{lem}[thm]{Lemma}
\newtheorem{prop}[thm]{Proposition}

\newtheorem{rem}[thm]{Remark}

\newtheorem*{acknowledgement}{Acknowledgements}

\numberwithin{figure}{section}

\newcommand{\bs}{\symbol{92}}
\newcommand{\pnc}{\mathbb{C} \mathbb{P}^n}
\newcommand{\plc}{\mathbb{C} \mathbb{P}^2}
\newcommand{\drc}{\mathbb{C} \mathbb{P}^1}
\newcommand{\plr}{\mathbb{R} \mathbb{P}^2}
\newcommand{\cc}{\mathbb{C}}
\newcommand{\rr}{\mathbb{R}}
\newcommand{\zz}{\mathbb{Z}}

\newcommand{\Tet}{\Theta}
\newcommand{\TetCh}{\Theta_{\Ch}}
\newcommand{\TetEP}{\Theta_{\EP}}
\newcommand{\TetChP}{\Theta_{\Ch}^P}
\newcommand{\TetCS}{\Theta_{\CS}}
\newcommand{\tTetCh}{\widetilde{\Theta}_{\Ch}}
\newcommand{\tTetChP}{\widetilde{\Theta}_{\Ch}^P}
\newcommand{\TetBC}{\Theta_{\BC}}

\newcommand{\ctg}{T^{\ast}}

\newcommand{\sph}{\mathbb{S}}
\newcommand{\cerc}{\mathbb{S}^1}
\newcommand{\SxS}{\mathbb{S}^2 \times \mathbb{S}^2}

\newcommand{\tGamma}{\widetilde{\Gamma}}
\newcommand{\tE}{\widetilde{E}}

\newcommand{\rEP}{\rho_{\EP}}
\newcommand{\rCh}{\rho_{\Ch}}
\newcommand{\rBC}{\rho_{\BC}}

\DeclareMathOperator{\im}{Im}

\DeclareMathOperator{\Id}{Id}

\DeclareMathOperator{\FS}{FS}
\DeclareMathOperator{\Ch}{Ch}
\DeclareMathOperator{\EP}{EP}
\DeclareMathOperator{\CS}{CS}
\DeclareMathOperator{\BC}{BC}

\newcommand{\FubSt}{\omega_{\FS}}

\title{On exotic monotone Lagrangian tori in $\plc$ and $\SxS$} 
\author{Agn\`es GADBLED}
\date{March 2011}

\begin{document}

\maketitle

\begin{abstract}
In this note, we prove that two constructions of exotic monotone 
Lagrangian tori, namely the one by Chekanov and Schlenk 
(see~\cite{Chek-SchI}, \cite{Chek-Sch-MSRI}) and the one obtained 
by the circle bundle construction of Biran 
(see~\cite{Bir-Cor-uniruling}), are Hamiltonian isotopic in~$\plc$ 
and $\SxS$.
\footnotetext{MSC classification: 53D05, 53D12, 57R17.}
\end{abstract}

\section*{Introduction}
\label{sec:intro}

Because of Darboux's theorem and because any open subset 
of~$\cc^n$ contains a Lagrangian torus, it is possible to construct 
a Lagrangian torus in any symplectic manifold.
The construction of Lagrangian submanifolds can become far more 
difficult as soon as we require some conditions - we will be here 
interested in monotonicity - on the Lagrangian submanifold.

Recall that a Lagrangian submanifold $L$ of a symplectic 
manifold~$(W,\omega)$ is said to be monotone (see Oh~\cite{Oh_I}) if 
there exists a non-negative constant~$K_L$ such that for every disc
$u \in \pi_2(W,L)$, 
$$\int u^\star \omega = K_L \; \mu_L (u),$$
where $\mu_L (u)$ is the Maslov class of the disc $u$.
Recall also that the existence of a monotone Lagrangian submanifold 
implies that the ambient symplectic manifold itself must be monotone,
which means that there exists a non-negative constant $K_W$ such 
that for every sphere $v~\in~\pi_2(W)$, 
$$\int v^\star \omega = K_W \; c_1(v),$$
where $c_1$ is the first Chern class of $(W,\omega)$.
Note that because of the relation between the first Chern class and 
the Maslov class of a sphere in~$\pi_2(W)$, in the case the first 
Chern class of $W$ does not vanish identically on $\pi_2(W)$, we 
have the following relation between the monotonicity constants:
$$K_W = 2 K_L.$$
In particular in this case, the monotonicity constant of any Lagrangian 
submanifold is prescribed by the monotonicity constant of the ambient 
symplectic manifold, and this gives already many restrictions on the 
Lagrangian.

Even in the case of $\cc^n$ endowed with its canonical symplectic 
structure, not many constructions of monotone Lagrangian tori are 
available.
The first and simplest example is the Clifford (or split) torus, 
that is the 
product of $n$~circles enclosing disks of the same area.
As any construction of a Lagrangian torus gives rise to infinitely many 
tori with the same symplectic invariants by applying a Hamiltonian 
diffeomorphism, we are rather interested in equivalence classes of tori 
under the action of the group of Hamiltonian diffeomorphisms.
It was only in 1995 that Chekanov~(\cite{Chekanov96}) gave the first 
examples of monotone Lagrangian tori in $\cc^n$ which are not 
Hamiltonian isotopic to the Clifford torus. 
Such tori are said exotic.

The Clifford torus can be embedded in the complex projective space  
and the product of spheres via the embedding of a ball or a polydisc of 
suitable size in~$\pnc$ or in~$\times_n \sph^2$. It is a monotone 
Lagrangian torus in $\pnc$ or in~$\times_n \sph^2$ called again 
Clifford torus. The Clifford torus was also the only known example 
of monotone Lagrangian torus in~$\pnc$ and in~$\times_n \sph^2$ 
till Chekanov and Schlenk~(\cite{Chek-SchI}, 
\cite{Chek-Sch-MSRI}) described in~2006 their families of 
exotic monotone Lagrangian tori.

Biran et Cornea (\cite{Bir-Cor-uniruling}) have also given recently a 
construction of monotone Lagrangian tori thanks to the circle bundle 
construction of Biran~(\cite{Biran-NonIntersections}).
It is believed that one should recover the families of monotone 
Lagrangian tori of Chekanov and Schlenk with the constructions 
through circle bundles. The present note proves 
this fact for $\plc$ and $\SxS$ and we will deal with the higher 
dimensions in some future work.

The tori of Chekanov and Schlenk in $\plc$ and $\SxS$ are defined by 
modifying slightly the construction of the exotic monotone Lagrangian 
torus in $\cc^2$ of Eliashberg and Polterovich in~\cite{El-Pol-1993}:
given a ``suitable'' (see Section~\ref{sec:CS-torus} 
and~\ref{sec:CS-torus-insphere}) curve 
$\gamma : [0, 2 \pi] \longrightarrow \cc$, 
consider 
$$\left\{{ \left. \left({\gamma(s) e^{i\theta}, \gamma(s) e^{-i\theta} 
}\right) \right| \theta \in [0, 2 \pi], s \in [0, 2 \pi] } \right\}$$
and then embed it into $\plc$ or $\SxS$.

The monotone torus in $\plc$ (or analogously in~$\SxS$) defined by 
Biran and Cornea in~\cite{Bir-Cor-uniruling} is not coming from a 
torus in $\cc^2$ and an embedding of a ball or a polydisc: one 
starts with a symplectic sphere $\Sigma$ (coming from a polarisation 
of $\plc$ or $\SxS$, see~\cite{Biran-Barriers}), and constructs 
the torus as the restriction to an equator of the sphere $\Sigma$ 
of a circle subbundle of a disc bundle~$E_{\Sigma}$ over $\Sigma$.
In the case of $\SxS$, this construction is known to be equal to the 
first example in the litterature of exotic monotone Lagrangian torus 
given by Entov and Polterovich in~\cite{En-Pol-Rigid}.

In this note we prove the following:

\begin{thmintro}
The construction of monotone exotic torus of Chekanov and Schlenk and 
the one given by the circle bundle construction of Biran are 
Hamiltonian isotopic in $\plc$ and $\SxS$.
\end{thmintro}

The idea of the proof for the projective plane is to relate 
Biran and Cornea's torus to the 
first construction of exotic monotone torus in $\cc^2$ given by 
Chekanov in~\cite{Chekanov96}.
This construction in $\cc^2$ is known to be Hamiltonian isotopic to 
the monotone exotic torus of Eliashberg and Polterovich 
in~\cite{El-Pol-1993}.
A similar isotopy between this tori can be used in the case of~$\plc$ 
to define an isotopy between Chekanov and Schlenk's torus and a 
modified Chekanov's torus.
After embedding into~$\plc$, one can prove that this modified 
Chekanov's torus is a circle subbundle of the disc bundle~$E_{\Sigma}$
 over an equator of $\Sigma$ and then isotope this torus inside the 
disc bundle to the torus of Biran and Cornea. 
The strategy for~$\SxS$ is analogous.

The article is organised the following way: in the first section we 
recall the two constructions of exotic monotone Lagrangian torus in 
$\cc^2$, namely the first one by Chekanov~\cite{Chekanov96} and the one 
of Eliashberg and Polterovich~\cite{El-Pol-1993}, and we exhibit a 
Hamiltonian isotopy between the two.
In the second section, we focus on the case of~$\plc$, recall the 
constructions of Chekanov and Schlenk~(\cite{Chek-SchI}, 
\cite{Chek-Sch-MSRI}), and Biran and Cornea~(\cite{Bir-Cor-uniruling}), 
and prove the existence of the Hamiltonian isotopy in this case. The 
third section deals with the case of~$\SxS$.

\begin{acknowledgement}
\emph{This work follows from a question raised by Felix Schlenk. I wish 
to thank him, Paul Biran, Octav Cornea and Leonid Polterovich for their 
explanations and remarks on these constructions. I also want to thank 
Kai Cieliebak, Ivan Smith and Maksim Maydanskyi for helpful 
discussions.}
\end{acknowledgement}

\noindent
This work was carried out while the author held a postdoctoral 
fellowship attached to the SNF project "Complexity and recurrence 
in Hamiltonian systems" at the University of Neuch\^atel, a 
Mathematical Sciences Research Institute postdoctoral fellowship, 
and a postdoctoral fellowship at the University of Cambridge 
funded by the European Research Council 
grant ERC-2007-StG-205349. 
The author wishes to thank these three institutions for their 
hospitality and stimulating research atmosphere.

\section{The first constructions in $\cc^2$}
\label{sec:C2}

In this section, we recall the first constructions of monotone 
exotic tori in~$\cc^2$ as they will be useful to understand the 
constructions in the compact symplectic manifolds $\plc$ and 
$\SxS$.

In order to simplify the computations and to avoid too many constants, 
we will use normalisations of the standard symplectic forms 
on~$T^\ast M$ for $M = \rr^{n}$ or $M = \cerc$ (the circle $\cerc$ 
being identified with $\rr / 2 \pi \zz$) such that the Liouville 
$1$-form on $T^\ast M$ is: 
\begin{equation}
\label{normalisation}
\lambda = \dfrac{1}{\pi} \Sigma p_i d q_i,
\end{equation}
where $(p,q)$ are the usual local coordinates on cotangent bundles, 
$q=(q_1, \ldots, q_n)$ coordinates on the basis and 
$p=(p_1, \ldots, p_n)$ coordinates in the fibers.
For example, with such a normalisation, the integral of the 
Liouville form along the circle centered at the origin and of 
radius~$r$ in $T^\ast \rr \simeq \rr^2$ is $r^2$.

\subsection{The first description by Chekanov}

Chekanov has given a first description of exotic monotone tori in
$\rr^{2n}$ in~\cite{Chekanov96}. We recall here his construction.

For any Lagrangian submanifold $L$ in $\rr^{2n}$, he has defined a 
Lagrangian submanifold $\Tet (L)$ in $\rr^{2n+2}$
the following way.
Consider the embedding 
$$
\begin{array}{cccc}
i_n: & \cerc \times \rr^n & \longrightarrow & \rr^{n+1} \\
     & (\theta,x_1,\ldots,x_n) & \longmapsto     & (e^{x_1}
\cos(\theta),e^{x_1} \sin(\theta),x_2,\ldots,x_n).
\end{array}
$$
Then $I_n = (i_n^\ast)^{-1}$ is a symplectic embedding of 
$T^\ast (\cerc \times \rr^n)$ into $T^\ast \rr^{n+1}$.
Denote by $\Tet (L)$ the image by the symplectic embedding~$I_n$ 
of the product of the zero section $N_0$ of $T^\ast S^1$ with 
the Lagrangian submanifold $L$. If $L$ is monotone in $\rr^{2n}$, 
then~$\Tet(L)$ is monotone in $\rr^{2n+2}$ with the same 
monotonicity constant.

We are here particularly interested in the case $n=1$. 
We give the expression in
coordinates of the map $I = I_1$ as it will be useful in the
following. 
If $\theta \in \cerc, \tau \in T^\ast_{\theta} \cerc, x \in \rr$ and $y \in
T^\ast_x \rr$, 
then 
$$I ((\tau,t) , (y,x)) = (p_0, p_1, q_0, q_1)$$ 
with 
$(q_0, q_1) \in \rr^2$ and $(p_0, p_1) \in T_{(q_0, q_1)} \rr^2$ such
that:
$$ \left \{ { 
\begin{array}{ccc}
q_0 & = & e^x \cos (\theta) \\
q_1 & = & e^x \sin (\theta) \\
p_0 & = & e^{-x} ( - \tau \sin (\theta) + y \cos(\theta) ) \\
p_1 & = & e^{-x} ( \tau \cos (\theta) + y \sin(\theta) ).
\end{array}
} \right.
$$
Identifying $T^\ast \rr^2$ with $\cc^2$ via the symplectomorphism 
$$
\begin{array}{ccc}
T^\ast \rr^2         & \longrightarrow & \cc^2 \\
(q_0, q_1, p_0, p_1) & \longmapsto     & (q_0 + i p_0, q_1 + i p_1),
\end{array}
$$
the map $I$ can be written as the following embedding of $T^\ast(\cerc
\times \rr)$ into $\cc^2$: 
$$\left({(\tau,\theta) , (y,x)} \right) \mapsto 
\left ( { 
\begin{array}{ccc}
z_0 & = & (e^{x} + i e^{-x} y) \cos(\theta) - i \tau e^{-x} \sin(\theta)\\
z_1 & = & (e^{x} + i e^{-x} y) \sin(\theta) + i \tau e^{-x} \cos(\theta)
\end{array}
} \right).
$$

For $n=1$, Chekanov's construction with $L$ the circle centered at the
origin of area $2 r^2$ can be parametrized in $\cc^2$ by:
\begin{equation}
\label{eq:ParamChek}
 \left \{ { 
\begin{array}{ccc}
z_0 & = & (e^{x} + i e^{-x} y) \cos(\theta)\\
z_1 & = & (e^{x} + i e^{-x} y) \sin(\theta)
\end{array}
} \right.
\end{equation}
with $(x,y) \in L$. 
It is a monotone Lagrangian torus with monotonicity constant $r^2$.
In \cite[Theorem 4.2]{Chekanov96}, Chekanov proved that this torus 
is not Hamiltonian isotopic to the Clifford torus (by versal 
deformations, using Ekeland-Hofer capacities and the displeacement 
energy). In the following, we will call it Chekanov's torus and denote 
it~$\TetCh (r^2)$.

\subsection{The version by Eliashberg and Polterovich and its relation 
with the previous construction}

Eliashberg and Polterovich have given in~\cite{El-Pol-1993} another
description of an exotic monotone torus in $\rr^4 \simeq \cc^2$.

If $D$ is a disk of $\cc$ of area $r^2$, $r>0$, which does 
not contain the origin, and $c = \{c(s) \; | \; s \in [0, 2 \pi]\}$ is its
boundary, parametrised by a smooth map 
$c~:~[0,2\pi]~\rightarrow~\cc$, then the torus defined as:
$$\left\{{ \left. \left({c(s) e^{i\theta}, c(s) e^{-i\theta} }\right) 
\right| \theta \in [0, 2 \pi], s \in [0, 2 \pi] } \right\}$$
is a monotone torus of monotonicity constant $r^2$ denoted 
$\TetEP(r^2)$. Eliashberg and Polterovich have proved that this torus 
is again not Hamiltonian isotopic to a split Lagrangian torus 
(\cite[Proposition 4.2.B]{El-Pol-1993}) by counts of holomorphic 
discs with boundary along this torus.

\begin{prop} 
\label{prop:EPetCh}
For any positive radius $r$, the monotone torus $\TetEP(r^2)$ is 
Hamiltonian isotopic to $\TetCh(r^2)$ in $\cc^2$.
\end{prop}

We give here a detailed proof of this well known result as we 
will use it in the next sections. We fix a positive radius $r$ and to 
simplify the notations, we will drop $r$ from the notations in the 
proof.

This proposition can be deduced from a series of lemmata.\\

The exotic torus of Eliashberg and Polterovich is constructed as the
orbit of some circle under a Hamiltonian circle action, so that it
satisfies the following:

\begin{lem}
\label{action EP}
The torus $\TetEP$ is stable under the following
action $\rEP$ of the circle on $\cc^2$: for $\theta \in [0;2 \pi]$ and
$(z_0,z_1) \in \cc^2$, 
$$\rEP(e^{i \theta})(z_0,z_1) = 
\left ( { \begin{array}{cc} e^{i \theta} & 0 \\ 0 & e^{-i \theta}
    \end{array}  } \right ) \left ( { \begin{array}{c} z_0 \\ z_1
    \end{array}  } \right )
=  \left ( { \begin{array}{c} e^{i \theta}  z_0 \\ e^{-i \theta} z_1
    \end{array}  } \right )
$$
of Hamiltonian 
$$H(z_0,z_1) = \dfrac{1}{2 \pi} \left({ |z_0|^2-|z_1|^2 } \right).$$
\end{lem}

More precisely, the torus $\TetEP$ is the orbit of the curve 
$$C = \left \{ \left .{ \left ( { \begin{array}{c} c(s) \\ c(s)
    \end{array}  } \right ) } \right | s \in  [0, 2 \pi] \right \}$$
under the action $\rEP$.\\

The exotic torus $\TetCh$ satisfies a similar property:

\begin{lem}
\label{lem:actionCh}
The torus $\TetCh$ is stable under the following action $\rCh$ of the
circle on $\cc^2$: for $\theta \in [0;2 \pi]$ and $(z_0,z_1) \in
\cc^2$,  
$$\rCh(e^{i \theta})(z_0,z_1) = 
\left ( { \begin{array}{cc}  \cos(\theta) & - \sin(\theta) \\
      \sin(\theta) & \cos(\theta) \end{array}  } \right ) \left ( {
    \begin{array}{c} z_0 \\ z_1  \end{array}  } \right )
$$
of Hamiltonian 
$$H(z_0,z_1) = \dfrac{1}{\pi} \im \left({z_0 \bar{z}_1}\right).$$
\end{lem}

In particular, the parametrisation~(\ref{eq:ParamChek}) gives that the 
torus $\TetCh$ is the orbit under the action $\rCh$
of the following curve of $\cc^2$:
$$\left \{ { \left . { \left ( { \begin{array}{c} e^x+ i e^{-x} y \\ 0
    \end{array}  } \right ) } \right | (x,y) \in L } \right \}$$
where $L$ is the circle of $\cc$ centered in the origin of
radius~$r$. But this torus can also be described as the
orbit of the curve:
$$\Lambda = \left \{ { \left . { \left ( { \begin{array}{c}
            \frac{1}{\sqrt{2}}(e^x+ i e^{-x} y) \\
            - \frac{1}{\sqrt{2}}(e^x+ i e^{-x} y) 
    \end{array}  } \right ) } \right | (x,y) \in L } \right \}.$$
\ \\

\begin{lem}
\label{lem:conjug} The two Hamiltonian actions $\rEP$ and $\rCh$ are
conjugate inside the special unitary group $SU(2)$.
\end{lem}

\begin{proof}
If we denote by $P$ the matrix
$$P = 
\left ( { \begin{array}{cc}  \frac{1}{\sqrt{2}}i & - \frac{1}{\sqrt{2}} \\
      \frac{1}{\sqrt{2}} & - \frac{1}{\sqrt{2}}i \end{array}  }
\right )$$
then $P$ is a matrix of the special unitary group such that for every
$\theta \in [0; 2 \pi]$,
$$\,^t \bar{P} 
\left ( { \begin{array}{cc} e^{i \theta} & 0 \\ 0 & e^{-i \theta}
    \end{array}  } \right )  P
= 
\left ( { \begin{array}{cc}  \cos(\theta) & - \sin(\theta) \\
      \sin(\theta) & \cos(\theta) \end{array}  } \right ). $$
\end{proof}

Now let $\TetChP$ be the torus obtained as the orbit under the
action $\rEP$ of the curve $P \Lambda$ of $\cc^2$.
Then $\TetChP$ is the image by the map $P$ of the torus $\TetCh$.

Note moreover that the diffeomorphism $P$ of $\cc^2$ is a Hamiltonian
diffeomorphism (because $P$ is a matrix of $SU(2)$). This means that 
$\TetChP$ is Hamiltonian isotopic to $\TetCh$ and it is now sufficient
for Proposition~\ref{prop:EPetCh} to prove that $\TetChP$ is
Hamiltonian isotopic to $\TetEP$. \\

\begin{lem}
\label{lem:ChPetEP}
The tori $\TetChP$ and $\TetEP$ are Hamiltonian isotopic inside $\cc^2$.
\end{lem}

To prove this, we will use the following lemma:

\begin{lem}[See \cite{Polt-Geom-2001}]
\label{lem:lemmePolterovich} 
  A Lagrangian isotopy is exact if and only if it can be extended to
  an ambient Hamiltonian isotopy.
\end{lem}

Recall (\cite[Section 6.1]{Polt-Geom-2001}) that given a Lagrangian 
isotopy of a closed manifold $N$ into a symplectic manifold 
$(W,\omega)$:
$$\Phi : N \times [0;1] \rightarrow W,$$
the pull-back $\Phi^\star \omega$ is of the form $\alpha_s \wedge ds$
where $\{\alpha_s\}$ is a family of closed $1$-forms on $N$. The
Lagrangian isotopy is said to be exact if $\alpha_s$ is exact for
all~$s$. \\

\noindent \emph{Proof of Lemma~\ref{lem:ChPetEP}.} \;
Thanks to Lemma~\ref{lem:lemmePolterovich}, it is enough to prove
that the tori $\TetChP$ and $\TetEP$ are exact Lagrangian isotopic. 

We note first that 
$$\Lambda^P = P \Lambda =  
 \left \{ { \left . \frac{1}{\sqrt{2}}(e^x+ i e^{-x} y)
\frac{1}{\sqrt{2}} (1 + i)
{ \left ( { \begin{array}{c}
            1 \\
            1 
    \end{array}  } \right ) } \right | (x,y) \in L } \right \},$$
so that $C$ and $\Lambda^P$ are both curves lying in the
diagonal of $\cc^2$.

Moreover, as $c$ is the boundary of the disc $D$, the integral of 
the Liouville form of $\cc$ is equal to $r^2$ on the curve $c$.

On the other hand, let $f$ be the map 
$$
\begin{array}{cccc}
f: & T^\ast \rr & \longrightarrow & \cc \\
   & (x,y)      & \longmapsto     & e^{x}+ i e^{-x} y.
\end{array}
$$
The map $f$ is an exact symplectomorphism (i.e. it preserves the 
Liouville form $p dq$) from $T^\ast \rr$ onto its image 
$$\{z \in \cc \; | \; \Re e(z) > 0 \}.$$
In particular, the integral of the Liouville form on the curve $f(L)$ 
is equal to the integral of the Liouville form on $L$. The rotation 
$z~\mapsto~\frac{1}{\sqrt{2}} (1+i) z$ preserves the Liouville 
form of~$\cc$ so that the integral of the Liouville form is also equal 
to~$r^2$ on the curve 
$$L^P=\frac{1}{\sqrt{2}} f(L) \frac{1}{\sqrt{2}} (1 + i).$$

As $L^P$ and $c$ are two closed curves in $\cc$ on which the integral of 
the Liouville form takes the same value (that is they are the boundary of
domains of the same area), these two curves can be Hamiltonianly
isotoped one into the other in $\cc$ (this can also be seen by applying 
Lemma~\ref{lem:lemmePolterovich}).
We denote by $\varphi: \cc \times [0;1] \longrightarrow \cc$ a Hamiltonian
isotopy of $\cc$ such that:
$$\varphi_0 = \Id \mbox{ and } \varphi_1\left({L^P}\right) = c.$$
As $L^P$ and $c$ are boundary of domains which do not contain
the origin, the Hamiltonian isotopy can be chosen so that:
\begin{equation}
\label{eq:nonsing}
\forall \; t \in [0;1], \; \varphi_t(0) = 0,
\end{equation}
(in other words, $\varphi_t\left({L^P}\right)$ never 
crosses the origin).

The two curves $\Lambda^P$ and $C$ are then Hamiltonian isotopic 
inside the diagonal $\Delta_{\cc^2}$ of $\cc^2$ via the isotopy
$(\varphi, \varphi)$, and we use this Hamiltonian isotopy and the 
circle action $\rEP$ to construct an exact Lagrangian isotopy 
from~$\TetChP$ to~$\TetEP$. This Lagrangian isotopy 
$\Phi: \TetChP \times [0;1] \longrightarrow \cc^2$ is
defined for $t \in [0;1]$ by:
$$
\begin{array}{cccc}
\Phi_t: & \TetChP  & \longrightarrow & \cc^2 \\
      & \left ( { \begin{array}{c}
            e^{i \theta} z\\
            e^{-i \theta} z
    \end{array}  } \right ) &
\longmapsto & \left ( { \begin{array}{c}
            e^{i \theta} \varphi(z,t)\\
            e^{-i \theta} \varphi(z,t)
    \end{array}  } \right )
\end{array}
$$
and it is a well defined Lagrangian isotopy because of
property~(\ref{eq:nonsing}).
As $\varphi$ is a Hamiltonian isotopy and $\rEP$ is a 
Hamiltonian circle action, one can check that $\Phi$ is 
exact. Thanks to Lemma~\ref{lem:lemmePolterovich}, this 
ends the proof of Lemma~\ref{lem:ChPetEP} and proves 
Proposition~\ref{prop:EPetCh}. \hfill \qedsymbol

\section{In the complex projective plane}

In this section, the projective complex plane $\plc$ will be 
endowed with the Fubini-Study symplectic form normalized so 
that the area of any complex projective line is 
$$\int_{\drc} \FubSt = 2.$$
This implies in particular that:
\begin{itemize}
\item The monotonicity constant of $\plc$ is $\frac{2}{3}$ and 
  any monotone Lagrangian submanifold in $\plc$ has monotonicity 
  constant $\frac{1}{3}$.
\item With the normalisation (\ref{normalisation}) of the symplectic 
  form of $\rr^4 \simeq \cc^2$, the open ball $B(2)$ of radius
  $\sqrt{2}$ is symplectically embedded into $\plc$ via the canonical
  embedding:
  \begin{equation}
  \label{embball}
  \begin{array}{cccc}
E_2:&  B(2)     & \longrightarrow & \plc \\
    &(z_0,z_1)& \longmapsto     & 
\left[ \, {z_0:z_1:\sqrt{2-|z_0|^2-|z_1|^2} } \, \right]
  \end{array}
  \end{equation}
\end{itemize}

Because of the symplectic embedding~(\ref{embball}), a way to 
construct a monotone torus into $\plc$ is to begin with a monotone 
torus of monotonicity constant~$\frac{1}{3}$ in $B(2)$ and embed 
this torus into $\plc$ via $E_2$. This is for example a possible 
construction of the Clifford torus: one begins with the product of 
the two circles centered in the origin and of 
radius~$\sqrt{\frac{2}{3}}$ in each factor of~$\cc^2 = \cc \times \cc$.

One could try to embed also Chekanov's torus as a monotone Lagrangian 
torus of $\plc$ via $E_2$. In order to get the right monotonicity 
constant, one should begin with $L$ a circle of area $\frac{2}{3}$. 
Unfortunately, one can check that the torus $\TetCh(\frac{2}{3})$ 
does not sit in the open ball $B(2)$ (for example with the 
parametrisation~(\ref{eq:ParamChek}), for $s=0$ and any $\theta$, 
one has $\|z\| > \sqrt{2}$).

Eliashberg and Polterovich's torus $\TetEP$ as it has been defined in
Section~\ref{sec:C2} is also not lying in the ball $B(2)$, so that it
cannot be embedded inside the complex projective plane $\plc$. But this 
construction can be modified in order to define an exotic torus in $\plc$.

\subsection{Chekanov and Schlenk's torus in the complex projective
  plane}
\label{sec:CS-torus}

Instead of defining a torus as the orbit under $\rEP$ of the boundary 
of a disk which does not contain the origin, one can also define a 
torus as the orbit of any closed embedded curve which is the boundary 
of a domain which does not contain the origin of $\cc$. The first 
torus of the family of monotone tori in~$\pnc$ defined by Chekanov 
and Schlenk~(\cite{Chek-SchI}, \cite{Chek-Sch-MSRI}) is constructed this way.

\begin{figure}[htbp]
\label{fig:gamma}
  \begin{center}
   \psfrag{C}{$\cc$}
   \psfrag{1}{$1$}
   \psfrag{0}{$0$}
   \psfrag{gamma}{$\color{red} \mathbf{\gamma}$}
   \includegraphics{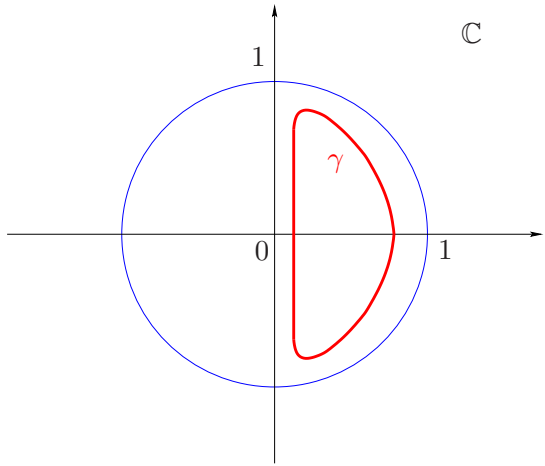}
   \caption{The curve $\gamma$} 
  \end{center}
\end{figure}

To be precise, replace the circle $c$ in the Eliashberg-Polterovich 
construction by the curve $\gamma$ which is the boundary
of a domain of area $\frac{1}{3}$ sitting inside the disk of
radius~$1$ centered in the origin of $\cc$ and inside the half-plane
of complex numbers of positive real part (see Figure \ref{fig:gamma}). 

Then the
curve 
$$\Gamma = \left ( {\begin{array}{c}
            \gamma\\
            \gamma
    \end{array}   } \right )$$
is lying inside the ball $B(2)$ and as the action of the cercle $\rEP$
is by multiplication by matrices of $SU(2)$, the orbit of this curve
is also entirely contained in $B(2)$. The torus constructed this way
can then be embedded in~$\plc$ through the canonical 
embedding~(\ref{embball}) and will be denoted~$\TetCS$. Chekanov and 
Schlenk have proved (see \cite{Chek-SchII}, \cite{Chek-Sch-MSRI})) that 
this torus is not Hamiltonian isotopic to the Clifford torus in~$\plc$ 
and is not Hamiltonianly displaceable in $\plc$.

\subsection{The torus of Biran and Cornea}

The Biran-Cornea torus is defined using the Lagrangian circle bundle
construction of Biran (see~\cite{Biran-Ciel-Subcrit}
and~\cite{Biran-NonIntersections}).

Let $\Sigma$ be the quadric of $\plc$, given by the homogeneous
equation:
$$z_0^2+z_1^2+z_2^2=0.$$
Then $(\plc, \omega, J ; \Sigma)$ (where $J$ is the standard complex
structure on $\plc$) is a polarized K\"ahler manifold of degree $1$ in
the sense of~\cite{Biran-Barriers} and~\cite{Biran-NonIntersections},
which means that $\Sigma$ is a smooth and reduced complex hypersurface
of the K\"ahler manifold $(W, \omega, J)$ (with $W = \plc$) whose 
homology class $[\Sigma]~\in~H_2(W,\zz)$ represents the Poincar\'e 
dual of~$[\omega]$.

In the case of such a polarisation, Biran proved
in~\cite{Biran-Barriers} that there exists an isotropic $CW$-complex
$\Delta~\in~(W,\omega)$ whose complement $(W \bs \Delta, \omega)$ 
is symplectomorphic to a standard symplectic disc bundle
$(E,\omega_0)$ modeled on the normal bundle $N_\Sigma$ of $\Sigma$ in
$W$ and whose fibres have area $1$.

The standard symplectic disc bundle is by definition the open unit
disc subbundle of the complex line bundle $\pi: N_\Sigma \rightarrow
\Sigma$ (with respect to some Hermitian metric), endowed with the
symplectic structure $\omega_0$ given by the following formula:
$$\omega_0 = \pi^\star \omega_{|\Sigma} + d(r^2 \alpha),$$
where $r$ is the radial coordinate on the fibres of $E_\Sigma
\rightarrow \Sigma$ defined using the Hermitian metric and $\alpha$ is
a connection $1$-form on $E_\Sigma \bs \Sigma$ with curvature
$$d \alpha = - \pi^\star(\omega_\Sigma).$$

In the case $\Sigma$ is the quadric above, $\Delta$ is $\mathbb{R}
\mathbb{P}^2$.\\

Biran and Cornea's torus is then constructed the following way
(see~\cite[Section 6.4.1]{Bir-Cor-uniruling}): the quadric $\Sigma$ is
topologically a sphere, so that we can consider an
equator~$\mathcal{E}$ of~$\Sigma$, that is a circle which divides the
sphere into two discs of equal areas. Let now $\TetBC$ be the
restriction to the equator $\mathcal{E}$ of the circle subbundle of
radius $\sqrt{\frac{2}{3}}$ inside the disc bundle $E_\Sigma
\rightarrow \Sigma$. This defines an exotic monotone torus in 
$\plc$ whose Floer homology (with $\zz/2$ coefficients) is trivial.

The following remark of Biran on the description of $E_\Sigma$ for the 
quadric will be crucial in the proof that $\TetBC$ is Hamiltonian 
isotopic to $\TetCS$ in~$\plc$. It was also used recently by Opshtein 
to prove a symplectic embedding Theorem (see~\cite[Theorem 4]{Ophstein}).

\begin{prop}
\label{prop:normal}
 Let $a, b,$ and $c$ be three real numbers, not identically zero, and 
let $D_{a,b,c}$ denote the projective line of $\plc$ defined by the 
homogeneous equation
\begin{equation}
\label{D}
 a z_0 + b z_1 + c z_2 = 0.
\end{equation}
The line $D_{a,b,c}$ intersects the quadric $\Sigma$ in two complex 
conjugate points $x$ and $\bar{x}$ and $\plr$ cuts this line along 
a circle in two complex conjugate discs which are the fibers of 
$E_\Sigma$ in $x$ and $\bar{x}$, respectively.
\end{prop}

\begin{proof}
The projective line $D_{a,b,c}$ is a curve of degree $1$ and 
intersects the curve $\Sigma$ of degree $2$ in two points. 
As the line and the quadric are defined by equations with real 
coefficients, their intersection points are complex conjugate. 
Moreover, they are distinct because the quadric has no real point.

Let $x = [x_0, x_1, x_2]$ be an intersection point of $D_{a,b,c}$ 
and $\Sigma$. The tangent space $T_x$ at the quadric in this point $x$ 
is defined by the homogeneous equation:
\begin{equation}
 \label{T}
x_0 z_0 + x_1 z_1 + x_2 z_2 = 0.
\end{equation}
In the complex linear space $\cc^3$, the two equations~(\ref{D}) 
and~(\ref{T}) define two complex hyperplanes $\tilde{D}$ and 
$\tilde{T}$ respectively. They are distinct as $(x_0, x_1, x_2)$ 
is not a multiple of a real vector. As a consequence, they have 
an intersection of complex dimension~$1$.
But the vector $(x_0, x_1, x_2)$ of $\cc^3$ is a non-zero vector 
in the intersection of $\tilde{D}$ and $\tilde{T}$ so that 
$$\tilde{D} \cap \tilde{T} = \mbox{Vect}_\cc (x_0, x_1, x_2).$$
Moreover the vector $(\bar{x}_0, \bar{x}_1, \bar{x}_2)$ belongs 
to $\tilde{D}$ and is orthogonal to $(x_0, x_1, x_2)$ so that 
$\tilde{D}$ is the direct orthogonal sum of the two vector 
spaces spanned by $(x_0, x_1, x_2)$ and 
$(\bar{x}_0, \bar{x}_1, \bar{x}_2)$, respectively. 
Note also that the vector $(\bar{x}_0, \bar{x}_1, \bar{x}_2)$ is 
in the Hermitian-ortogonal to the vector space $\tilde{T}$. 
Consequently, the complex projective lines $D_{a,b,c}$ and $T_x$ are 
orthogonal for the induced Hermitian product on $\plc$.

The complex projective line $D_{a,b,c}$ is stable under complex 
conjugation and its invariant submanifold $D_{a,b,c} \cap \plr$ 
is a circle which divides $D_{a,b,c}$ in two discs which are 
images of each other by the complex conjugation, and thus have 
the same area (equal to $1$ with our conventions).  Moreover, 
each disc contains only one of the intersection points with the 
quadric $x$ or~$\bar{x}$.

To be sure that the two discs above are exactly the fibers 
of~$E_\Sigma$, one should be sure that the discs are ``centered'' 
in the intersection points.
One way to prove this is to see that it is true in the case $a=1$, 
$b=0$ and $c=0$ and it can be extended to the other cases thanks 
to a transformation of~$U(3)$. \end{proof}

\begin{rem}
Through any point of~$\cc^3$ passes a complex line 
defined by an equation with real coefficients (
for dimensional reason, the Hermitian orthogonal supplement 
of a point must intersect $\rr^3$). As a consequence, 
Proposition~\ref{prop:normal} describes the bundle $E_\Sigma$ 
in any point of the quadric $\Sigma$.
\end{rem}

\subsection{Biran and Cornea's torus is Hamiltonian isotopic to Chekanov and Schlenk's torus}

\begin{thm}
 The Biran and Cornea's torus is Hamitonian isotopic to Chekanov and 
Schlenk's torus in $\plc$.
\end{thm}

\begin{proof}
First we also define a modified Chekanov torus using the 
curve~$\gamma$. Let~$\tGamma$ be the curve 
$$\tGamma = \left ( {\begin{array}{c}
            \gamma\\
            - \gamma
    \end{array} } \right )$$
and $\tTetCh$ the torus defined as the orbit of the curve $\tGamma$
under the action of~$\rCh$.
Note that as the curve $\gamma$ is Hamiltonian isotopic to the 
curve~$L$ in $\cc$, the same proof as the one of 
Lemma~\ref{lem:ChPetEP} proves
that $\tTetCh$ is Hamiltonian isotopic to $\TetCh$ in $\cc^2$.

The same way as for $\TetCS$, because $\tGamma$ is contained in $B(2)$ 
and because the action of the cercle $\rCh$ is by multiplication by 
matrices of $SU(2)$, the torus $\tTetCh$ is entirely contained in $B(2)$
and can be embedded in $\plc$ as a monotone torus.

Let also $\tTetChP$ be the torus defined as the orbit under the action
$\rEP$ of the curve $P \Gamma$. As before, 
$$\tTetChP = P \tTetCh,$$
so that $\tTetChP$ and $\tTetCh$ are Hamiltonian isotopic. Moreover,
the Hamiltonian isotopy is given by a family of matrices of $U(2)$ so
that this isotopy is preserving the ball $B(2)$. 
As a consequence, this isotopy enables to define an isotopy between 
$\tTetChP$ and $\tTetCh$ inside the complex projective plane.

Now notice that in $B(2)$,
$$\TetCS = R_{-\frac{\pi}{4}} \tTetChP$$
where $R_{-\frac{\pi}{4}}$ is the rotation $e^{- i \frac{\pi}{4}}$ in 
each factor of $\cc^2$, so that $\TetCS$ is the image of $\tTetChP$ 
by a Hamiltonian diffeomorphism. 

To prove the theorem, we will prove that $\tTetCh$ is Hamiltonian 
isotopic to~$\TetBC$.

The torus~$\tTetCh$ can be defined in homogeneous coordinates 
of~$\plc$ as the set:
$$\tTetCh = \left \{ \left. [\cos(\theta) \sqrt{2} \gamma(s): \sin(\theta) \sqrt{2} \gamma(s) : \sqrt{2 - 2 |\gamma(s)|^2}] \; \right| \; \theta, s \in [0, 2 \pi] \right\}.$$
It is Hamiltonian isotopic in $\plc$ to any torus obtained 
rotation inside the plane of the torus~$\tTetCh$, and in 
particular to the torus 
$$e^{i \frac{\pi}{2}} \tTetCh$$
described in $\plc$ as the set
$$\tTetCh' = \left \{ \left. [\cos(\theta) \sqrt{2} \gamma(s): \sin(\theta) \sqrt{2} \gamma(s) : i \sqrt{2 - 2 |\gamma(s)|^2}] \; \right| \; \theta, s \in [0, 2 \pi] \right\}.$$
The latter can also be seen as the image of $\tTetCh$ by the 
following modified embedding of the ball $B(2)$ inside $\plc$:
 \begin{equation}
  \label{embball2}
  \begin{array}{cccc}
\tE_2:&  B(2)     & \longrightarrow & \plc \\
      & (z_0,z_1)& \longmapsto     & 
\left[ \, {z_0:z_1: i \sqrt{2-|z_0|^2-|z_1|^2} } \, \right].
  \end{array}
  \end{equation}
It is thus enough to prove that the torus 
$\TetBC$ is Hamiltonian isotopic to $\tTetCh'$ in $\plc$ and 
for that purpose, it will be more convenient to work with the 
embedding~$\tE_2$.

Let $Z$ be the cylinder of $T^\ast \cerc \times \{0\} \subset \ T^\ast (\cerc \times \rr)$ defined as:
$$Z = \left \{ {(\theta,\tau,0,0) \in T^\ast (\cerc \times \rr^n) \;|\; |\tau| < 1 } \right \}.$$
The image of this cylinder by Chekanov's map $I$ can be parametrized 
in~$T^\ast \rr^2$ by:
$$ \left \{ { 
\begin{array}{ccc}
q_0 & = & \cos (\theta) \\
q_1 & = & \sin (\theta) \\
p_0 & = & - \tau \sin (\theta) \\
p_1 & = & \tau \cos (\theta) 
\end{array}
} \right.
$$
for $|\tau| < 1$ and in $\cc^2$ by:
$$\left ( { 
\begin{array}{ccc}
z_0 & = & \cos(\theta) - i \tau \sin(\theta)\\
z_1 & = & \sin(\theta) + i \tau \cos(\theta)
\end{array}
} \right).
$$
As $|\tau| < 1$, any point in the image by $I$ of the cylinder $Z$ is 
lying inside the ball $B(2)$. Its image $I(Z)$ under the embedding 
(\ref{embball2}) is parametrized by:
$$\{ [\cos(\theta) - i \tau \sin(\theta): \sin(\theta) + i \tau \cos(\theta) : i \sqrt{1-\tau^2}] \; | \; s \in [0, 2 \pi], \tau \in (-1,1) \}$$

\begin{lem}
 The quadric $\Sigma$ is the union of the image $I(Z)$ and the two points 
``at infinity'' $[1:i:0]$ and $[1:-i:0]$.
\end{lem}

\begin{proof}
Any point of $I(Z)$ satisfies the equation of the quadric. Moreover, when 
$\tau$ tends to $\varepsilon \in \{-1;1\}$, the point of homogeneous 
coordinates 
$$[\cos(\theta) - i \tau \sin(\theta): \sin(\theta) + i \tau \cos(\theta) : 
i \sqrt{1-\tau^2}]$$ 
tends to the point $[1 : \varepsilon : 0]$, so that $I(Z)$ can be 
compactified in the closed surface $\Sigma$.
\end{proof}

\begin{lem}
 The image by $I$ of the zero section $\tau = 0$ is an equator of the 
 quadric. 
\end{lem}

Indeed, the zero section $\tau = 0$ of $T^\ast \cerc$ cuts the cylinder 
$Z$ in two pieces of equal area so that its image cuts the quadric 
in two discs of same area $1$, one disc containing the point 
$[1 : i : 0]$ and the other one the point $[1 : -i : 0]$. In the 
homogeneous coordinates of $\plc$, this equator can be parametrised by:
$$[\cos(\theta):\sin(\theta):i],\; \theta \in [0 ; 2 \pi].$$

\begin{lem}
 For a fixed $\theta \in \rr$, the curve 
$$\{[\cos(\theta) \sqrt{2} \gamma(s): \sin(\theta) \sqrt{2} \gamma(s) : i \sqrt{2 - 2 |\gamma(s)|^2}]\; | \; s \in [0, 2 \pi] \}$$
lies inside the fiber of $E_\Sigma$ at the point 
$[\cos(\theta):\sin(\theta):i]$.\\
Moreover, this curve enclose a domain of area $\frac{2}{3}$ inside the 
fiber at the base point $[\cos(\theta):\sin(\theta):i]$.
\end{lem}

\begin{proof}
Thanks to Proposition~\ref{prop:normal}, we know how to describe the 
disc bundle~$E_\Sigma$ at a point of the quadric as soon as we have 
the equation of a real projective line containing this point. 
In the case of a point $[\cos(\theta):\sin(\theta):i]$ on the equator 
of the $\Sigma$, it lies on the real line of homogeneous equation
$$\sin(\theta) z_0 - \cos(\theta) z_1 = 0.$$
Because of Proposition~\ref{prop:normal}, this means that the 
fibers of the disc bundle $E_\Sigma$ in the points
$$[\cos(\theta):\sin(\theta):i] \mbox{ and }
[\cos(\theta):\sin(\theta):-i]=[\cos(\theta+\pi):\sin(\theta+\pi):i]$$
are the two discs of the projective line $N_\theta$ of equation
 $$\sin(\theta) z_0 - \cos(\theta) z_1 = 0$$
separated by the curve $N_\theta \cap \plr$.

Any point of the curve 
$$\{[\cos(\theta) \sqrt{2} \gamma(s): \sin(\theta) \sqrt{2} \gamma(s) : i \sqrt{2 - 2 |\gamma(s)|^2}] \; | \; s \in [0, 2 \pi] \}$$
belongs to the line $N_\theta$.
Moreover, this curve does not intersect $\plr$ as the curve $\gamma$ 
does not intersect neither the imaginary axes nor the circle of radius~$1$. 
Therefore, this curve is entirely contained in one of the fiber of 
the disc bundle $E_\Sigma$.
Actually, the image of the half disc
$$\{z \in \cc \; | \; \Re e (z) > 0 \mbox{ and } |z|< 1\}$$
by the map
$$z \mapsto \{[\cos(\theta) \sqrt{2} z: \sin(\theta) \sqrt{2} z : i \sqrt{2 - 2 |z|^2}]\}$$
is contained in $N_\theta$ and does not intersect $\plr$. Moreover the
area of the image of the half disc is equal to $1$ (that is twice the 
area of the original half disc) so that it is the entire fiber of $E_\Sigma$ 
at the point $[\cos(\theta):\sin(\theta):i]$, image of the complex number 
$\frac{1}{\sqrt{2}}$.
The curve $\gamma$ is enclosing a domain of area~$\frac{1}{3}$, and 
consequently the curve 
$$\{[\cos(\theta) \sqrt{2} \gamma(s): \sin(\theta) \sqrt{2} 
\gamma(s) : i \sqrt{2 - 2 |\gamma(s)|^2}] \; | \; s \in [0, 2 \pi] \}$$
is enclosing a domain of area $\frac{2}{3}$.
\end{proof}

To end the proof, one can see that in any fiber of the normal bundle, 
the curve 
$$\{[\cos(\theta) \sqrt{2} \gamma(s): \sin(\theta) \sqrt{2} \gamma(s) : i \sqrt{2 - 2 |\gamma(s)|^2}] \; | \; s \in [0, 2 \pi] \}$$
is isotopic to the circle of area $\frac{2}{3}$ (used 
to construct the torus of Biran and Cornea). In order to construct 
a Lagrangian isotopy between $\tTetCh'$ and $\TetBC$, one can use 
the isotopy defined for the fiber of the point corresponding to 
$\theta = 0$ and then extend this isotopy by the action of the 
circle~$\rCh$. Lemma~\ref{lem:lemmePolterovich} enables then to 
extend this exact Lagrangian isotopy into a Hamiltonian isotopy 
of~$\plc$.
\end{proof}

\section{In the product of two two-dimensional spheres}

Let us now consider $W = \SxS = \drc \times \drc$ 
endowed with the symplectic form 
$$\omega_W = \FubSt \oplus \FubSt$$
where $\FubSt$ is normalized such that 
$$\int_{\drc} \FubSt = 1,$$
which means that $\drc$ is symplectomorphic to the sphere of 
radius~$\frac{1}{2}$ in~$\rr^3$ with our conventions.
With this normalisation, the product of the two spheres is monotone 
with monotonicity constant 
$$K_W = \frac{1}{2}.$$
As a consequence, the monotonicity constant of any monotone
Lagrangian submanifold in this product is 
$$K_L = \frac{1}{4}.$$

In this symplectic manifold, one has also the two corresponding constructions 
of exotic monotone Lagrangian torus.

\subsection{Chekanov and Schlenk's torus in the product of spheres}
\label{sec:CS-torus-insphere}

For the construction by Chekanov and Schlenk, one has this time to 
begin with a curve $\gamma$ enclosing a domain of 
area~$\frac{1}{4}$ in the complex 
plane inside the half-disc of radius~$1$ and positive real part.

With this choice of curve $\gamma$ the torus of $\cc^2$ obtained 
by action of $\rEP$
 $$\left\{{ \left. \left({\gamma(s) e^{i\theta}, \gamma(s)
 e^{-i\theta} }\right) \right| \theta \in [0, 2 \pi], 
s \in [0, 2 \pi] } \right\}$$
is contained into the product $B(1) \times B(1)$, where $B(1)$ is 
the ball of radius~$1$ in $\cc$.
As this product can be symplectically embedded into 
the product $\drc \times \drc$ via 
 \begin{equation}
  \label{embprodball}
  \begin{array}{cccc}
E_{1,1}: &  B(1) \times  B(1) & \longrightarrow & \drc \times \drc \\
        & (z_0,z_1)         & \longmapsto     & 
\left({ \left[ \, {z_0:\sqrt{1-|z_0|^2}} \, \right] , \left[ \, {z_1:\sqrt{1-|z_1|^2}} \, \right] } \right),
  \end{array}
  \end{equation}
this construction produces a monotone Lagrangian torus
still denoted~$\TetCS$ in~$\drc \times \drc$ which is not Hamiltonian 
isotopic to the Clifford torus and 
non displaceable (see \cite{Chek-Sch-MSRI}, \cite{Chek-SchII}).

\subsection{The torus constructed by Biran's circle bundle construction}

The construction of Biran involves again a polarisation. This 
time, the symplectic hypersurface $\Sigma$ is the diagonal 
sphere described as
$$\Sigma = \{ (x,x) \in  \SxS | \; x \in \sph^2\}$$
in $\SxS$, or as the hypersurface of $\drc \times \drc$ satisfying 
the homogeneous equation:
$$z_0 w_1 = z_1 w_0,$$
where $[z_0:w_0]$ are the homogeneous coordinates on the first 
copy of $\drc$ and $[z_1:w_1]$ on the second.

The total space of the standard symplectic disc 
bundle~$(E_\Sigma,\omega_0)$ modeled on the normal 
bundle~$N_\Sigma$ of $\Sigma$ in $W$ and whose fibres have 
area $1$ is in the case of $\drc \times \drc$ symplectomorphic 
to the complement of the antidiagonal:
$$\Delta = \{ (x,-x) \in  \SxS | \; x \in \sph^2\}$$
described in $\drc \times \drc$ as the Lagrangian submanifold 
of homogeneous equation:
$$z_0 \bar{z}_1 + w_0 \bar{w}_1 = 0.$$
The exotic torus, still denoted $\TetBC$, is now defined as the 
restriction of the circle subbundle of radius $\frac{1}{\sqrt{2}}$ 
over an equator of $\Sigma$.

This torus is known to be equal to the exotic monotone Lagrangian 
torus of Entov and Polterovich (\cite[Exemple 1.22]{En-Pol-Rigid}, 
see also the study in~\cite[Section~2]{El-Polt-QS}):
$$K = \left \{ {\left(x,y) \right) \in \SxS \; \left|{ \; x_3 + y_3 = 0, x_1 y_1 + x_2 y_2 + x_3 y_3 = - \frac{1}{2} }\right. } \right \},$$
with $\sph^2$ being in their conventions the sphere of radius $1$ of 
$\rr^3$, the symplectic form being rescaled.

It is also equal to the torus defined in the cotangent bundle of 
$\sph^2$ by the geodesic flow (see~\cite{MR2417887}) and embedded 
in a suitable way in a Weinstein's neighbourhood of the Lagrangian 
sphere $\Delta$.

\subsection{The two constructions are Hamiltonian isotopic}

\begin{thm}
Biran's exotic torus is Hamitonian isotopic to Chekanov and 
Schlenk's exotic torus in $\SxS$.
\end{thm}

\begin{proof}
In the case of the product of spheres, we cannot use the 
corresponding modified Chekanov's torus as in the case of the 
complex projective plane because this torus does not embed into 
the product of balls of radius $1$.

However, in $\cc^2$, the modified Chekanov's torus $\tTetCh$ sits 
in the normal bundle of the image $I(Z)$ of the cylinder 
$Z$ of the 
original description, and more precisely in the fibers over the 
image of the zero section. 
We also know from Section~\ref{sec:CS-torus} 
that $\TetCS$ is the image of $\tTetCh$ by the composition of the 
rotations $R_{-\frac{\pi}{4}}$ and the Hamiltonian diffeomorphism 
described by the matrix~$P$. Therefore it also sits in the normal 
bundle of a surface $Q$, namely 
the image of $I(Z)$ by $R_{-\frac{\pi}{4}} P$.

This surface $Q$ can be parametrized, with the 
coordinates $(\tau, \theta)$ coming from the original 
coordinates on $\ctg \cerc$ by:
$$\left ( { 
\begin{array}{ccc}
z_0 & = & \left(\frac{1}{\sqrt{2}} - \frac{1}{\sqrt{2}} \tau \right) e^{i (\theta + \pi/4)}\\
z_1 & = & \left(\frac{1}{\sqrt{2}} + \frac{1}{\sqrt{2}} \tau \right) e^{-i (\theta + \pi/4)}
\end{array}
} \right).
$$
This surface is, as expected, invariant by the action of 
$\rEP$ and the torus~$\TetCS$ (for the moment considered 
in~$\cc^2$) lies in the normal bundle of the 
surface along the curve $\tau = 0$. Moreover, 
as $\tau \in (0,1)$, $Q$ is contained in the 
product~$B(1) \times  B(1)$ and can be embedded into 
$\SxS$.

We would like to relate this surface to the diagonal $\Sigma$ as 
in the proof of the $\plc$ case. However, the diagonal $\Sigma$ 
is not invariant by the action of the cercle $\rEP$ but by the following action of the circle which can be written on~$\cc^2$:
$$\rBC(e^{i \theta})(z_0,z_1) = 
\left ( { \begin{array}{cc} e^{i \theta} & 0 \\ 0 & e^{i \theta}
    \end{array}  } \right ) \left ( { \begin{array}{c} z_0 \\ z_1
    \end{array}  } \right )
$$
This action cannot be conjugate into $SU(2)$, but it will be enough 
to conjugate the two circle actions inside $SO(3) \times SO(3)$ using 
the description of $W$ as the product of two spheres of $\rr^3$.
As the rotation $e^{i \theta}$ on $\drc$ corresponds to the 
rotation of matrix
$$ \left ( { \begin{array}{ccc} 
\cos \theta & -\sin \theta & 0 \\
\sin \theta &  \cos \theta & 0 \\
        0   &       0      & 1
    \end{array}  } \right )$$
one sees that $\rEP$ and $\rBC$ are conjugate in 
$SO(3) \times SO(3)$, for example by the pair $(P_1,P_2)$
where $P_1$ is the identity matrix and 
$$ P_2 = \left ( { \begin{array}{ccc} 
1 & 0 & 0 \\
0 &-1 & 0 \\
0 & 0 & -1
\end{array}  } \right ).$$

In order to see the action of the diffeomorphism of matrix $P_2$ 
on $\sph^2$ we need to embed symplectically the ball $B(1)$ in the 
sphere of radius $\frac{1}{2}$ of $\rr^3$. Such a map can be given 
for example by
$$
 \begin{array}{cccc}
E_1:&  B(1)    & \longrightarrow & S^2 \\
    &  z=a+ib  & \longmapsto     & 
\left({ \begin{array}{ccc} \sqrt{1-|z|^2} \; a\\
                           \sqrt{1-|z|^2} \; b\\
                           \frac{1}{2} - |z|^2
 \end{array}} \right).
  \end{array}
$$
The image by $(P_1, P_2)$ of $(E_1,E_1)(Q)$ 
in $\SxS$ is now invariant by the action~$\rBC$, but is
not the diagonal $\Sigma$ (which was to be expected as the cylinder has an area equal to $4$):
$$\widetilde{Q}=(P_1,P_2)((E_1,E_1)(Q))= \left \{ (X,Y) \in \SxS \right \}$$
with $X$ parametrised by 
$$\left ( { 
\begin{array}{c}
\sqrt{\frac{1}{2}+\tau-\frac{1}{2} \tau^2} \left(\frac{1}{\sqrt{2}} - \frac{1}{\sqrt{2}} \tau \right) \cos (\theta + \pi/4)\\
\sqrt{\frac{1}{2}+\tau-\frac{1}{2} \tau^2} \left(\frac{1}{\sqrt{2}} - \frac{1}{\sqrt{2}} \tau \right) \sin (\theta + \pi/4)\\
\tau-\frac{1}{2} \tau^2
\end{array}
} \right)$$
and $Y$ parametrised by
$$\left ( { 
\begin{array}{c}
\sqrt{\frac{1}{2}-\tau-\frac{1}{2} \tau^2} \left(\frac{1}{\sqrt{2}} + \frac{1}{\sqrt{2}} \tau \right) \cos (\theta + \pi/4)\\
\sqrt{\frac{1}{2}-\tau-\frac{1}{2} \tau^2} \left(\frac{1}{\sqrt{2}} + \frac{1}{\sqrt{2}} \tau \right) \sin (\theta + \pi/4)\\
\tau + \frac{1}{2} \tau^2
\end{array}
} \right)
$$
for $s \in [0, 2 \pi]$ and $\tau \in (-1,1)$.

Now note that the image of the zero section (the curve $\tau = 0$) 
of the cylinder is an equator of the diagonal sphere. 
Moreover, the tangent bundle of~$\widetilde{Q}$ along 
the equator is equal to the tangent bundle of the sphere~$\Sigma$. 
This means that the original curve~$\gamma$ which was sitting 
in the normal bundle of~$Q$ is, after these Hamiltonian 
isotopies, a curve on which the integral of the Liouville form 
is equal to~$\frac{1}{2}$ sitting in the normal bundle of~$\Sigma$ 
along the equator. It is lying in the disc bundle of area one as it 
does not intersect the anti-diagonal~$\Delta$ after the Hamiltonian 
diffeomorphism.
One finishes the proof using Lemma~\ref{lem:lemmePolterovich} as in 
the case of~$\plc$.
\end{proof}

\bibliographystyle{plain}
\bibliography{Biblioanglaise}
\addcontentsline{toc}{section}{References}

\ \\
\begin{small}
Agn\`es GADBLED, 
Department of Pure Mathematics and Mathematical Statistics,
Centre for Mathematical Sciences,
University of Cambridge,
Wilberforce Road,
Cambridge,
CB3 0WB,
United Kingdom\\
electronic address: A.Gadbled@dpmms.cam.ac.uk
\end{small}

\end{document}